\documentclass{llncs}

\usepackage{color}
\usepackage{cite}
\usepackage{comment}
\usepackage{amsmath}
\usepackage{amsfonts}
\usepackage{graphicx} 
\usepackage{bm}
\usepackage{epstopdf}

\usepackage{xcolor,colortbl}
\usepackage{subfigure}
\allowdisplaybreaks

\definecolor{Gray}{gray}{0.85}

\begin{document}

\title{On principal component analysis of the convex combination of two data matrices and its application to acoustic metamaterial filters}

\author{Giorgio Gnecco$^1$, Andrea Bacigalupo$^2$}
\institute{$^1$ IMT School for Advanced Studies, Lucca, Italy \\
\email{giorgio.gnecco@imtlucca.it}\\
$^2$ University of Genoa, Italy \\
\email{andrea.bacigalupo@unige.it}
}


\maketitle

\begin{abstract}
In this short paper, a matrix perturbation bound on the eigenvalues 
found by principal component analysis is investigated, for the case in which the data matrix on which principal component analysis is performed is a convex combination of two data matrices. The application of the theoretical analysis to multi-objective optimization problems (e.g., those arising in the design of acoustic metamaterial filters) is briefly discussed, together with possible extensions. 
\end{abstract}


\section{Introduction}

Principal Component Analysis (PCA) is a well-known data dimensionality reduction technique \cite{Jolliffe2002}. It works by projecting a dataset of $m$ vectors $\mathbf{x}_j \in \mathbb{R}^n$, $j=1,\ldots,m$ (represented by a data matrix $\mathbf{X} \in \mathbb{R}^{m \times n}$, whose rows are such vectors) onto a reduced  $d$-dimensional subspace of $\mathbb{R}^n$, which is generated by the first $d<n$ so-called principal directions. These are orthonormal eigenvectors of the symmetric 
matrix $\mathbf{C} \doteq \frac{1}{m}\mathbf{X}' \mathbf{X} \in \mathbb{R}^{n \times n}$, associated with its $d$ largest positive eigenvalues. The latter are proportional (via the multiplicative factor $\frac{1}{m}$) to the $d$ largest positive eigenvalues of the related Gram matrix $\mathbf{G} \doteq \mathbf{X} \mathbf{X}' \in \mathbb{R}^{m \times m}$, whose element in position $(i,j)$ is the inner product between the vectors $\mathbf{x}_i$ and $\mathbf{x}_j$. Focusing on such eigenvalues is important because the eigenvalues corresponding to discarded principal directions (the ones associated with the successive eigenvalues, not selected by PCA) provide information about the mean squared error of approximation of the dataset when only the first $d$ principal directions are kept to construct that approximation (as a consequence, knowing these eigenvalues is useful also to select a suitable value for $d$). Moreover, when the dataset has zero mean, each such eigenvalue represents the empirical variance of the projection of the dataset onto the corresponding principal direction.

Given this framework, the goal of this short work is to get a matrix perturbation bound on the eigenvalues of the Gram matrix $\mathbf{G}$, for the case in which the data matrix $\mathbf{X}$ is a convex combination of other two data matrices. According to the authors' experience, this is a non-standard but potentially quite interesting way of using PCA. The application of the theoretical analysis to multi-objective optimization (a framework in which a convex combination of two data matrices can arise) is discussed in the last section, together with possible extensions.  




\section{Theoretical analysis}\label{sec:theoretical_analysis} 


In this section, a matrix perturbation bound on the eigenvalues of the Gram matrix of a dataset is provided, for the case in which its data matrix $\mathbf{X}(\alpha)$ is a convex combination of two data matrices $\mathbf{X}_1$ and $\mathbf{X}_2$, with varying weights $\alpha \in [0,1]$ and $1-\alpha$.\\   

{\bf Proposition.} {\em Let $\mathbf{X}_1, \mathbf{X}_2 \in \mathbb{R}^{m\times n}$ be two data matrices, $\alpha \in [0,1]$, $\mathbf{G}(\alpha) \doteq \mathbf{X}(\alpha)\mathbf{X}'(\alpha) \in \mathbb{R}^{m \times m}$ be the Gram matrix of their convex combination $\mathbf{X}(\alpha) \doteq \alpha \mathbf{X}_1 + (1-\alpha) \mathbf{X}_2$ with weights $\alpha$ and $1-\alpha$, $K \in \mathbb{N}$ and, for $k = 0,1,\ldots,K$, $\alpha_k \doteq \frac{k}{K}$. Let the non-negative eigenvalues of $\mathbf{G}(\alpha)$ and $\mathbf{G}(\alpha_k)$ be ordered,  
 respectively, as $\lambda_1(\mathbf{G}(\alpha))\geq \lambda_2(\mathbf{G}(\alpha)) \geq \ldots \geq \lambda_m(\mathbf{G}(\alpha))$ and $\lambda_1(\mathbf{G}(\alpha_k))\geq \lambda_2(\mathbf{G}(\alpha_k)) \geq \ldots \geq \lambda_m(\mathbf{G}(\alpha_k))$. Finally, let $\sigma_1(\mathbf{X}(\alpha))$ and $\sigma_1(\mathbf{X}(\alpha_k))$ be the largest singular values of $\mathbf{X}(\alpha)$ and $\mathbf{X}(\alpha_k)$, respectively. Then, for any $k=0,1,\ldots,K-1$ and $\alpha \in \left[\alpha_k,\alpha_{k+1} \right]$, the following holds, for all $i=1,\ldots,m$:
\begin{equation}\label{eq:bound}
|\lambda_i(\mathbf{G}(\alpha))-\lambda_i(\mathbf{G}(\alpha_k))|\leq \frac{2}{K} \left(\sigma_1(\mathbf{X}_1)+\sigma_1(\mathbf{X}_2)\right)^2\,.
\end{equation}
}
\begin{proof}
Using the singular value decomposition of $\mathbf{X}(\alpha)=\mathbf{U}(\alpha)\mathbf{\Sigma}(\alpha)\mathbf{V}'(\alpha)$ (being $\mathbf{U}(\alpha) \in \mathbb{R}^{m \times m}$ and $\mathbf{V}(\alpha) \in \mathbb{R}^{n \times n}$ orthogonal matrices, and $\mathbf{\Sigma}(\alpha) \in \mathbb{R}^{m \times n}$ a rectangular matrix whose $q \doteq \min\{m,n\}$ elements on its main diagonal are the singular values $\sigma_i(\mathbf{X}(\alpha))$, ordered from the largest singular value
to the smallest one), one gets
\begin{eqnarray}
\mathbf{G}(\alpha)&=&\mathbf{X}(\alpha)\mathbf{X}'(\alpha)=\mathbf{U}(\alpha)\mathbf{\Sigma}(\alpha)\mathbf{V}'(\alpha)\left(\mathbf{U}(\alpha)\mathbf{\Sigma}(\alpha)\mathbf{V}'(\alpha)\right)' \nonumber \\
&=&\mathbf{U}(\alpha)\mathbf{\Sigma}(\alpha)\mathbf{V}'(\alpha) \mathbf{V}(\alpha)\mathbf{\Sigma}'(\alpha)\mathbf{U}'(\alpha)=\mathbf{U}(\alpha)\mathbf{\Lambda}(\alpha)\mathbf{U}'(\alpha)\,,
\end{eqnarray}
where, denoting by $\mathbf{I}_{n \times n} \in \mathbb{R}^{n \times n}$ the identity matrix, the property $\mathbf{V}'(\alpha) \mathbf{V}(\alpha)=\mathbf{I}_{n \times n}$ has been used, and $\mathbf{\Lambda}(\alpha) \doteq \mathbf{\Sigma}(\alpha)\mathbf{\Sigma}'(\alpha) \in \mathbb{R}^{m \times m}$ is a diagonal matrix whose elements on its main diagonal are the squares $\sigma_i^2(\mathbf{X}(\alpha))$ of the singular values $\sigma_i(\mathbf{X}(\alpha))$ of $\mathbf{X}(\alpha)$, plus $m-n$ additional zeros (if and only if $m>n$). The $\sigma_i^2(\mathbf{X}(\alpha))$ and the possible $m-n$ additional zeros are also the eigenvalues $\lambda_i(\mathbf{G}(\alpha))$ of $\mathbf{G}(\alpha)$, since this is a symmetric and positive semi-definite matrix. We consider first the case $m\leq n$, then the case $m>n$.

\renewcommand{\labelitemi}{\textbullet}
\begin{itemize}
\item Case 1: $m \leq n$. We exploit the matrix perturbation bound on singular values provided in \cite[Theorem 3.3.16 (c)]{HornJohnson1990}, according to which, given any two matrices $\mathbf{A},\mathbf{B} \in \mathbb{R}^{m \times n}$, one has, for all $i=1,\ldots,\min\{m,n\}=m$,
\begin{equation}\label{eq:HornJohnson}
|\sigma_i(\mathbf{A})-\sigma_i(\mathbf{A}+\mathbf{B})| \leq \sigma_1(\mathbf{B})\,.
\end{equation}
Denoting by $\Delta\alpha$ a variation of $\alpha$, we apply Eq. (\ref{eq:HornJohnson}) with $\mathbf{A}=\mathbf{X}(\alpha)$ and $\mathbf{B}=\mathbf{X}(\alpha+\Delta\alpha)-\mathbf{X}(\alpha)=\Delta\alpha \left(\mathbf{X}_1 -\mathbf{X}_2\right)$. Recalling the relation $\lambda_i(\mathbf{G}(\alpha))=\sigma_i^2(\mathbf{X}(\alpha))$ valid for all $i=1,\ldots,m$, and the fact that $\sigma_i(\mathbf{X}(\alpha)) \geq 0$, we get
\begin{eqnarray}
&&|\lambda_i(\mathbf{G}(\alpha+\Delta\alpha))-\lambda_i(\mathbf{G}(\alpha))| \nonumber \\
&=&|\sigma_i^2(\mathbf{X}(\alpha+\Delta\alpha))-\sigma_i^2(\mathbf{X}(\alpha))| \nonumber \\
&=&|\sigma_i(\mathbf{X}(\alpha+\Delta\alpha))-\sigma_i(\mathbf{X}(\alpha))|\,(\sigma_i(\mathbf{X}(\alpha+\Delta\alpha))+\sigma_i(\mathbf{X}(\alpha))) \nonumber \\
&\leq& |\Delta\alpha| \,\sigma_1(\mathbf{X}_1-\mathbf{X}_2)\, (\sigma_i(\mathbf{X}(\alpha+\Delta\alpha))+\sigma_i(\mathbf{X}(\alpha)))\,.
\end{eqnarray}
Moreover, $0 \leq \sigma_1(\mathbf{X}_1-\mathbf{X}_2) \leq \sigma_1(\mathbf{X}_1)+\sigma_1(\mathbf{X}_2)$ by Eq. (\ref{eq:HornJohnson}) with $\mathbf{A}=\mathbf{X}_1$ and $\mathbf{B}=-\mathbf{X}_2$, whereas $0 \leq \sigma_i(\alpha \mathbf{X}_1+(1-\alpha)\mathbf{X}_2) \leq \alpha \sigma_i(\mathbf{X}_1)+(1-\alpha) \sigma_1(\mathbf{X}_2)\leq \sigma_1(\mathbf{X}_1)+\sigma_1(\mathbf{X}_2)$, where the second last inequality is obtained again by Eq. (\ref{eq:HornJohnson}) with $\mathbf{A}=\alpha \mathbf{X}_1$ and $\mathbf{B}=(1-\alpha)\mathbf{X}_2$. Combining all the above, one gets
\begin{equation}
|\lambda_i(\mathbf{G}(\alpha+\Delta\alpha))-\lambda_i(\mathbf{G}(\alpha))| \leq 2 |\Delta\alpha| \left(\sigma_1(\mathbf{X}_1)+\sigma_1(\mathbf{X}_2)\right)^2\,,
\end{equation}
from which one obtains Eq. (\ref{eq:bound}) for $|\Delta\alpha|=|\alpha_k-\alpha| \leq \frac{1}{K}$.

\item Case 2: $m > n$. The proof is the same as above for all but the last (smallest) $m-n$ eigenvalues of $\mathbf{G}(\alpha)$ and $\mathbf{G}(\alpha_k)$. However, the latter eigenvalues are all equal to $0$, and the bound (\ref{eq:bound}) still holds trivially for them.  
\end{itemize}
\end{proof}

The bound expressed by Eq. (\ref{eq:bound}), whose proof shows the Lipschitz continuity of the eigenvalues of $\mathbf{G}(\alpha)$ with respect to $\alpha$, can be used in the following way. First, one finds the sets of eigenvalues of the matrices $\mathbf{G}(\alpha_k)$, for $k=1,\ldots,K$. Then, for each $\alpha \in [0,1]$, one finds its nearest $\alpha_k$, then applies Eq. (\ref{eq:bound}) to locate approximately the eigenvalues of the new matrix $\mathbf{G}(\alpha)$.


\section{Discussion and possible extensions}

The theoretical framework considered in Section \ref{sec:theoretical_analysis} has  application, e.g., in the combination of PCA with the so-called weighted sum method, which is used in the context of multi-objective optimization \cite{ColletteSiarry2003}. In the case of two objective functions, this method approximates the Pareto frontier of a multi-objective optimization problem by minimizing, for $\mathbf{p} \in P \subseteq \mathbb{R}^n$, the trade-off
\begin{equation}
J_\alpha(\mathbf{p}) \doteq \alpha J_1(\mathbf{p}) + (1-\alpha) J_2(\mathbf{p})
\end{equation}
between the two objective functions $J_1(\mathbf{p})$ and $J_2(\mathbf{p})$, for different values of the parameter $\alpha \in [0,1]$ (an adaptive version of the method can be applied in cases for which the classical weighted sum method fails, e.g., when the Pareto frontier is nonconvex \cite{KimdeWeck2005}). Assuming that both $J_1(\mathbf{p})$ and $J_2(\mathbf{p})$ are differentiable and the optimization problem is unconstrained (i.e., $P=\mathbb{R}^n$) or that it can be reduced to an unconstrained optimization problem by using a suitable penalization approach, one could perform the optimization numerically by applying the classical gradient method, possibly combined with a multi-start approach. In order to reduce the computational effort needed for the exact computation of the gradient at each iteration of the gradient method, one could replace it with its approximation obtained by applying PCA to the gradient field $\nabla J_\alpha(\mathbf{p})$ evaluated on a subset of points $\mathbf{p}_j \in P$ (for $j=1,\ldots,m$), then projecting the exact gradient onto the subspace generated by the average of the gradients $\nabla J_\alpha(\mathbf{p}_j)$, and by the first principal directions found by PCA, when this is applied to the dataset $\{\nabla J_\alpha(\mathbf{p}_j)\}_{j=1}^m$, after a pre-processing step, which makes it centered\footnote{It is common practice to apply PCA to centered (also called de-meaned) data matrices $\mathbf{X}^{(c)}$, i.e., having the form $\mathbf{X}^{(c)} \doteq \mathbf{X}-\mathbf{1}_m \bar{\mathbf{x}}'$, where $\mathbf{1}_m \in \mathbf{R}^m$ denotes a column vector made of $m$ ones, and $\bar{\mathbf{x}}\in \mathbf{R}^n$ is a column vector whose elements are the averages of the corresponding columns of $\mathbf{X}$. This does not change the quality of the results of the theoretical analysis, because, by linearity, the centered convex combination of two data matrices $\mathbf{X}_1$ and $\mathbf{X}_2$ is equal to the convex combination of the two respective centered data matrices $\mathbf{X}_1^{(c)}$ and $\mathbf{X}_2^{(c)}$.}. Due to the structure of the objective function $J_\alpha(\mathbf{p})$, such dataset (represented by a data matrix $\mathbf{X}_\alpha$) would be made of the convex combination (with coefficients $\alpha$ and $1-\alpha$) of the two datasets $\{\nabla J_1(\mathbf{p}_j)\}_{j=1}^m$ and $\{\nabla J_2(\mathbf{p}_j)\}_{j=1}^m$, represented respectively by the two data matrices $\mathbf{X}_1$ and $\mathbf{X}_2$. In this context, the results of our theoretical analysis could be useful to restrict the application of the weighted sum method to a coarse grid of values $\alpha_k$ for $\alpha \in [0,1]$, from which one could infer, for other values of $\alpha$, the empirical variances of the projections of the (de-meaned) data matrices $\mathbf{X}^{(c)}(\alpha)$ onto the principal directions either selected or discarded by PCA, when PCA is applied to each such data matrix $\mathbf{X}^{(c)}(\alpha)$. Moreover, in view of this application to multi-objective optimization, the theoretical analysis of this work could be extended by finding upper bounds on the Jordan canonical angles\footnote{These, loosely speaking, represent the {\em smallest angles} between corresponding elements of the orthonormal bases of two subspaces of $\mathbb{R}^n$, being the bases chosen to minimize such angles. For rigorous definitions, see \cite{Wedin1972,ZhuKnyazev2013} and the references therein.} 
between the subspaces found by PCA applied to the data matrices $\mathbf{X}(\alpha)$ generated from $\mathbf{X}_1$ and $\mathbf{X}_2$ for two different values of $\alpha \in [0,1]$. Such an extension could be derived 
by applying a variation (proved in \cite{Wedin1972}) of the well-known Davis-Kahan theorem in matrix perturbation theory \cite[Theorem 3.4]{StewartSun1998}. 
A second extension of the analysis to the case of nonlinear versions of PCA, such as kernel PCA \cite{COAP2009}, seems also possible (e.g., via the kernel trick).

We conclude mentioning that, in our related work \cite{METANANO2021} about the design of acoustic metamaterial filters according to a single-objective optimization framework\footnote{Such optimization problems are typically characterized by a high computational effort needed for an exact evaluation of the gradient of their objective functions, which is motivated by the fact that each such evaluation requires solving the physical-mathematical model associated with the specific choice of the vector of parameters of the model, which is also the vector of optimization variables. The reader is referred to \cite{JOTA2020} for a further discussion about these computational issues.} (see \cite{Volterra2016} for a physical-mathematical model similar to the one considered in \cite{METANANO2021}), we have successfully applied PCA to the sampled gradient field of the objective function, achieving numerical results comparable with those obtained by using the exact gradient, but with a much smaller computational effort (e.g., with a reduction of the dimension by a factor $4$). A similar outcome is expected when moving to a multi-objective optimization framework. So, for this kind of optimization problems, the application of PCA to the approximation of the sampled gradient field of the objective function can be a valid alternative to the use of surrogate optimization methods (which replace the original objective function with a surrogate function, learned either offline \cite{JOTA2020} or online \cite{CMAME2021}), in case a gradient-based optimization algorithm is used to solve the optimization problem.

\vspace{-2mm}
\section*{Acknowledgment}
A. Bacigalupo and G. Gnecco are members of INdAM. The authors acknowledge financial support from INdAM-GNAMPA (project Trade-off between Number of Examples and Precision in Variations of the Fixed-Effects Panel Data Model), from INdAM-GNFM, from the Universit\`{a} Italo Francese (projects GALILEO 2019 no. G19-48 and GALILEO 2021 no. G21$\_$89), from the Compagnia di San Paolo (project MINIERA no. I34I20000380007), and from the University of Trento (project UNMASKED 2020).

\end{document}